\documentclass[11pt,a4paper]{amsart}

\usepackage[margin=1in]{geometry}
\usepackage[T1]{fontenc}
\usepackage{lmodern}
\usepackage{amsmath,amssymb,amsthm,mathtools}
\usepackage[protrusion=true,expansion=false]{microtype}
\usepackage{hyperref}
\hypersetup{colorlinks=true,linkcolor=blue,citecolor=blue,urlcolor=blue}

\newtheorem{theorem}{Theorem}
\newtheorem{lemma}{Lemma}

\newcommand{\C}{\mathbb C}
\newcommand{\D}{\mathbb D}
\newcommand{\Lem}{\Lambda_f}

\title[A degree-four lemniscate path theorem]{A degree-four lemniscate path theorem}
\author{Venkata Siddharth Pendyala}
\email{venkatasiddharthpendyala@gmail.com}

\begin{document}

\begin{abstract}
We prove the degree-four case of a path problem of Erd\H{o}s, Herzog, and Piranian.  If $f$ is monic of degree four and all zeros of $f$, counted with multiplicity, lie in the open unit disk, then two zeros from this list can be joined inside $\{z:|f(z)|<1\}$ by a possibly degenerate polygonal path of length less than $2$.
\end{abstract}

\maketitle

\section{The theorem}

Write
\[
        \D=\{z\in\C: |z|<1\},\qquad
        \overline{\D}=\{z\in\C: |z|\le1\},
\]
and, for $q\in\C$ and $R\ge0$,
\[
        \overline D(q,R)=\{z\in\C: |z-q|\le R\}.
\]
For a polynomial $f$, let
\[
        \Lem=\{z\in\C: |f(z)|<1\}.
\]
Throughout, zeros are listed with multiplicity.  Erd\H{o}s, Herzog, and Piranian asked whether, for every monic polynomial whose zeros lie in $\D$, the lemniscate sublevel set $\Lem$ contains a path of length less than $2$ joining two zeros of the polynomial~\cite{EHP}; see also Mac Lane~\cite{MacLane} and the modern problem record~\cite{Bloom1041}.  The following gives an elementary proof in degree four.

\begin{theorem}\label{thm:quartic}
Let
\[
        f(z)=\prod_{k=1}^4(z-z_k),\qquad z_1,z_2,z_3,z_4\in\D,
\]
where the zeros are listed with multiplicity.  Then there exist distinct indices $i,j$ such that $z_i$ and $z_j$ can be joined in $\Lem$ by a possibly degenerate polygonal path of length less than $2$.
\end{theorem}

The proof uses the following finite geometric lemma.

\begin{lemma}[Four-point radial lemma]\label{lem:fourpoint}
Let $w_1,w_2,w_3,w_4\in\overline{\D}$ satisfy
\[
        |w_i-w_j|\ge1\qquad(i\ne j).
\]
Then there are distinct indices $i,j$ such that
\[
        \prod_{k=1}^4 |tw_i-w_k|\le1\qquad(0\le t\le1)
\]
and
\[
        \prod_{k=1}^4 |tw_j-w_k|\le1\qquad(0\le t\le1).
\]
\end{lemma}

\begin{proof}
First suppose that all four points are nonzero.  Their arguments are distinct, since two nonzero points of $\overline{\D}$ on a common ray have distance less than $1$.  Order the points cyclically, and let
\[
        g_1,g_2,g_3,g_4>0,
        \qquad g_1+g_2+g_3+g_4=2\pi,
\]
be the four angular gaps.  Call a point balanced if the sum of its two adjacent gaps is at most $\pi$.  The adjacent-gap sums at the first and third points add to $2\pi$, and the adjacent-gap sums at the second and fourth points also add to $2\pi$.  Hence at least one point in each opposite pair is balanced, so there are two distinct balanced points.

It remains to show that every balanced point has a safe radial arm.  Rotate so that the balanced point is the positive real number $r$, where $0<r\le1$.  Its angular neighbours have the form
\[
        u=\rho e^{ia},\qquad v=\sigma e^{-id},
\]
with $0<\rho,\sigma\le1$, $a,d>0$, and $a+d\le\pi$.  Let $y$ be the fourth point, and put $x=tr$, where $0\le t\le1$.

We first compare a point with the boundary point on the same ray.  If $w=\rho e^{i\theta}\in\overline{\D}$, $0<\rho\le1$, and $|r-w|\ge1$, then
\begin{equation}\label{eq:comparison}
        |x-w|\le |x-e^{i\theta}|\qquad(0\le x\le r).
\end{equation}
Indeed,
\[
        |x-e^{i\theta}|^2-|x-\rho e^{i\theta}|^2
        =(1-\rho)(1+\rho-2x\cos\theta).
\]
For $\cos\theta\le0$ this is nonnegative.  For $\cos\theta>0$, the separation condition gives
\[
        2r\rho\cos\theta\le r^2+\rho^2-1.
\]
Dividing by $\rho$ and using $x\le r$ and $r\le1$, we get
\[
        2x\cos\theta\le2r\cos\theta
        \le \frac{r^2+\rho^2-1}{\rho}
        \le \rho\le1+\rho.
\]
Thus $1+\rho-2x\cos\theta\ge0$, and \eqref{eq:comparison} follows.

Applying \eqref{eq:comparison} to $u$ and $v$, we obtain
\[
        |x-u|\le \sqrt{1+x^2-2x\cos a},\qquad
        |x-v|\le \sqrt{1+x^2-2x\cos d}.
\]
Since $a+d\le\pi$,
\[
        \cos a+
        \cos d=2\cos\frac{a+d}{2}\cos\frac{a-d}{2}\ge0.
\]
A direct expansion gives
\begin{align*}
        &(1+x^2)^2-(1+x^2-2x\cos a)(1+x^2-2x\cos d)\\
        &\qquad=2x(1+x^2)(\cos a+\cos d)-4x^2\cos a\cos d.
\end{align*}
The right-hand side is nonnegative if $\cos a\cos d\le0$.  If instead $\cos a,\cos d\ge0$, then
\[
        2\cos a\cos d\le\cos a+\cos d.
\]
In this case the same right-hand side is at least
\[
        2x(\cos a+\cos d)(1-x+x^2)\ge0.
\]
Therefore
\[
        |x-u|\,|x-v|\le1+x^2.
\]
The remaining two factors satisfy $|x-y|\le1+x$ and $|x-r|=r-x$.  Hence
\[
        \prod_{k=1}^4 |x-w_k|
        \le (r-x)(1+x)(1+x^2)
        \le (1-x)(1+x)(1+x^2)
        =1-x^4\le1.
\]
Thus every balanced point gives a safe radial arm.

It remains to consider the case in which one of the four points is the origin.  The separation hypothesis then rules out any second point at the origin, and the other three points lie on $\partial\D$, because their distances from the origin are at least $1$ and their moduli are at most $1$.  The radial arm attached to the origin gives product zero.  Choose any one of the other three points and rotate it to $1$.  If the two remaining boundary points are $e^{i\alpha}$ and $e^{i\beta}$, then the product along this radial arm is
\[
        t(1-t)|t-e^{i\alpha}|\,|t-e^{i\beta}|,
        \qquad 0\le t\le1.
\]
Therefore
\[
        t(1-t)|t-e^{i\alpha}|\,|t-e^{i\beta}|
        \le t(1-t)(1+t)^2\le1,
\]
since $t(1-t)\le1/4$ and $(1+t)^2\le4$.  The origin and the chosen boundary point therefore give two distinct safe radial arms.
\end{proof}

\begin{proof}[Proof of Theorem~\ref{thm:quartic}]
If $|z_i-z_j|<1$ for some pair $i\ne j$, then the segment $[z_i,z_j]$, possibly degenerate, gives the required path.  Indeed, every point $z$ of this segment lies in $\D$, and
\[
        |z-z_i|\,|z-z_j|\le \frac{|z_i-z_j|^2}{4}.
\]
Each of the two remaining factors is strictly smaller than $2$, since both corresponding zeros and $z$ lie in $\D$.  Hence
\[
        |f(z)|< \frac{|z_i-z_j|^2}{4}\cdot2\cdot2
        =|z_i-z_j|^2<1,
\]
and the segment has length less than $1$.

We may therefore assume that
\[
        |z_i-z_j|\ge1\qquad(i\ne j).
\]
Let $\overline D(q,R)$ be the smallest closed disk containing the four zeros.  Since the zeros lie in $\D$, if
\[
        r_0=\max_{1\le k\le4}|z_k|,
\]
then $r_0<1$, and $\overline D(0,r_0)$ contains all four zeros.  Thus $0<R\le r_0<1$.  Define
\[
        w_k=\frac{z_k-q}{R},\qquad k=1,2,3,4.
\]
Then $w_k\in\overline{\D}$ and
\[
        |w_i-w_j|=\frac{|z_i-z_j|}{R}\ge \frac{1}{R}>1\qquad(i\ne j).
\]
By Lemma~\ref{lem:fourpoint}, there are distinct indices $i,j$ such that both radial arms from $0$ to $w_i$ and from $0$ to $w_j$ are safe.  Consider the broken segment
\[
        z_i\longrightarrow q\longrightarrow z_j.
\]
On the first arm, a point has the form $z=q+Rt w_i$ with $0\le t\le1$, and hence
\[
        |f(z)|
        =\prod_{k=1}^4 |q+Rt w_i-z_k|
        =R^4\prod_{k=1}^4 |tw_i-w_k|
        \le R^4<1.
\]
The same calculation applies to the second arm.  Its length is
\[
        |z_i-q|+|z_j-q|=R(|w_i|+|w_j|)\le2R<2.
\]
The broken segment is therefore contained in $\Lem$ and joins the two listed zeros $z_i$ and $z_j$ with length less than $2$.
\end{proof}

\end{document}